\RequirePackage{fix-cm}

\documentclass[smallextended]{svjour3}       

\smartqed  
\usepackage{}
\usepackage{amsfonts}
\usepackage{mathrsfs}
\usepackage[latin1]{inputenc}
\usepackage{version,tabularx,multicol}
\usepackage{stmaryrd}
\usepackage{amsmath}
\usepackage{amssymb}
\usepackage{mathabx}

\newtheorem{lemma}{Lemma}[section]
\newtheorem{theorem}{Theorem}[section]
\newtheorem{corollary}{Corollary}[section]
\newtheorem{definition}{Definition}[section]
\newtheorem{proposition}{Proposition}[section]

\newcommand{\beqlb}{\begin{eqnarray}}
\newcommand{\eeqlb}{\end{eqnarray}}
\newcommand{\beqnn}{\begin{eqnarray*}}
\newcommand{\eeqnn}{\end{eqnarray*}}

\def\qed{\hfill$\Box$\medskip}

\def\<{\langle}\def\>{\rangle}

 \def\ar{\!&}

\begin{document}

\title{Central limit theorems for supercritical superprocesses with immigration\thanks{Supported by NSFC 11301020 grants.}
}

\titlerunning{CLT for superprocesses with immigration}        

\author{Li Wang}

\authorrunning{L Wang}

\institute{Li Wang \at
School of Sciences, Beijing University of
Chemical Technology,
Beijing 100029,
P.R. China. \\
\email{wangli@mail.buct.edu.cn}
}


\maketitle

\begin{abstract}
In this paper, we establish a central limit theorem for a large class of general supercritical
superprocesses with immigration with spatially dependent branching mechanisms
satisfying a second moment condition. This central limit theorem extends and generalizes the
results obtained by Ren, Song and Zhang \cite{RSZ15}.
We first give law of large numbers for supercritical superprocesses with
immigration since there is few convergence result on immigration superprocesses, then based on
these results, we establish the central limit theorem.


\end{abstract}

\section{Introduction}

\setcounter{equation}{0}

In recent years, there have been many papers on law of large numbers type convergence theorems
and on central limit theorem types convergence
results for branching Markov processes and superprocess, see for instance, \cite{AH76b,AM15,CRW08,CS07,EJ09,EH10,EW06,LYS13,M12,RSZ14OU,RSZ14,RSZ15,WL10}
and the references therein.  Especially in \cite{RSZ15}, Ren, Song and Zhang proved one
central limit theorem which generalizes and unifies all the central limit theorems of \cite{M12,RSZ14OU} and the
advantage of this central limit theorem is that it allows us to characterize the limiting Gaussian
field and reveals more independent structures of the limiting Gaussian field.

The focus of this paper is, to take the next step, on central limit theorems for superprocesses with immigration.
First, we will give law of large numbers on supercritical superprocesses with immigration since
there is few convergence result on immigration superprocesses, then based on the
these results, we will give the central limit theorems. The main tool of this paper is the excursion
measure of superprocess and immigration superprocess. As a continuation of \cite{RSZ15},
the underlying spatial process in this paper is the same as in \cite{RSZ15}, we will state
it in the next subsection for reader's convenience.

\subsection{Spatial process}

\setcounter{equation}{0}

E is a locally compact separable metric space and $m$ is a $\sigma$-finite Borel measure on $E$ with full
support. $\partial$ is a point not contained in $E$ and will be interpreted as the cemetery point. Every
function $f$ on $E$ is automatically extended to $E_\partial :=E\cup\{\partial\}$ by setting $f(\partial)= 0$.
We will assume that $\xi=\{\xi_t,\Pi_x\}$ is an $m$-symmetric Hunt process on E and $\zeta:=\inf\{t > 0 :\xi_t=\partial\}$
is the lifetime of $\xi$. The semigroup of $\xi$ will be denoted by $\{P_t : t \geq 0\}$. We will always assume that
there exists a family of continuous strictly positive
symmetric functions $\{p_t(x,y), t>0\}$ on $E\times E$ such that
 \beqnn
P_tf(x)=\int_E p_t(x,y)f(y)m(dy).
 \eeqnn
It is well known that, for $p\geq 1$, $\{P_t:t\geq 0\}$ is a strongly continuous contraction semigroup on $L^p(E,m)$.
Define $\widetilde{a}_t(x):=p_t(x,x)$. We will always assume that
$\widetilde{a}_t(x)$ satisfies the following two conditions:
\begin{enumerate}
\item[(a)]For any $t>0$, we have
\beqnn
\int_E\widetilde{a}_t(x)m(dx)<\infty.
\eeqnn

\item[(b)]There exists $t_0>0$ such that for all $t\geq t_0$, $\widetilde{a}_t(x)\in L^2(E,m)$.

\end{enumerate}
These two conditions are satisfied by a lot of Markov processes, see the examples in \cite{RSZ14}.

\subsection{Superprocess with immigration}

Let $\mathcal{B}_b(E)$ ($\mathcal{B}^+_b(E)$) be the set of (positive) bounded Borel measurable
functions on $E$. Let $X=\{X_t:t\geq 0\}$ be a superprocess determined by the following
parameters: a spatial motion $\xi=\{\xi_t,\Pi_x\}$ satisfying the assumptions of the previous
subsection, a branching rate function $\beta(x)$ on $E$ which is a non-negative bounded measurable
function and a branching mechanism $\psi$ of the form
 \beqnn
\psi(x,\lambda)=-a(x)\lambda+b(x)\lambda^2+\int_{(0,+\infty)}(e^{-\lambda y}-1+\lambda y)n(x,dy),
~x\in E,~\lambda>0,
 \eeqnn
where $a\in \mathcal{B}_b(E)$, $b\in \mathcal{B}^+_b(E)$ and $n$ is a kernel from $E$ to $(0,\infty)$
satisfying
 \beqnn
\sup_{x \in E}\int_0^\infty y^2n(x,dy)<\infty.
 \eeqnn
Let $M(E)$ be the space of finite measures on $E$ equipped with the topology of
weak convergence. The existence of such superprocess is well known, see, for instance, \cite{D93} or \cite{Li11}.
$X$ is a c\'{a}dl\'{a}g Markov process taking values in $M(E)$. For any $\mu\in M(E)$,
we denote the law of $X$ with initial configuration $\mu$ by $\mathbb{P}_\mu$. Then $X$ has
transition semigroup $(Q_t)_{t\geq 0}$  defined by
\beqnn
\int_{M(E)}e^{-\nu(f)}Q_t(\mu,d\nu)=\exp\{-\mu(V_t(f))\},
\eeqnn
where $V_t(f)$ is the unique positive solution to the equation
\beqlb\label{1.1}
V_t(f)(x)+\Pi_x\int_0^t\psi(\xi_s,V_{t-s}f(\xi_s))\beta(\xi_s)ds=\Pi_xf(\xi_t),
\eeqlb
where $\psi(\partial,\lambda)=0, \lambda>0.$
Let $M(E)^\circ:=M(E)\setminus \{0\}$, where $0$ is the null measure.
Define
 \beqlb\label{1.2}
\varphi(f)=\eta(f)+\int_{M(E)^\circ}(1-e^{-\nu(f)})H(d\nu),~f\in \mathcal{B}^+_b(E),
 \eeqlb
where $\eta\in M(E)$ and $\nu(1)H(d\nu)$ is a finite measure on $M(E)^\circ$.

Suppose that $\{(Y_t, \mathscr{G}_t): t\geq 0\}$ is a Markov process
in $M(E)$ with transition semigroup
$(Q^N_t)_{t\geq 0}$ given by
\beqlb\label{1.3}
\int_{M(E)}e^{-\nu(f)}Q^N_t(\mu,d\nu)=\exp\left\{-\mu(V_t(f))-\int_0^t\varphi(V_s(f))ds\right\},
\eeqlb
where $V_t(f)$ is the unique positive solution to the equation (\ref{1.1}). We call $\{(Y_t, \mathscr{G}_t): t\geq 0\}$
a superprocess with immigration associated with $(Q_t)_{t\geq 0}$
or $\{X_t:t\geq 0\}$. For the exsitence
and properties of such processes, we refer the reader to \cite{Li92,LiS95,Li96}.
For any $\mu\in M(E)$,
we still use $\mathbb{P}_\mu$ to denote the law of $Y$ with initial configuration $\mu$ when there is
no confusion.

\medskip

Define the random set
\beqnn
\mathcal{Z}:= \overline{\{t\geq0: \|Y_t\|=0\}}.
\eeqnn
In this paper, we assume that the process $X$, or equivalently $Y$, is supercritical, then clearly  $\mathbb{P}_{\mu}\{\exists~t>0: \|Y_t\|=0\}=0$ or
$\mathcal{Z}$ is bounded $\mathbb{P}_\mu$-a.s. for $\mu\in M(E)$.

\medskip

Let
\beqnn
\alpha(x):=\beta(x)a(x)~~\mbox{and}~~A(x):=\beta(x)\left(2b(x)+\int_0^\infty y^2n(x,dy)\right).
\eeqnn
Then, by our assumptions, $\alpha(x)\in \mathcal{B}_b(E)$, $A(x)\in \mathcal{B}^+_b(E)$. Thus there exists $M>0$ such that
\beqnn
\sup_{x\in E}(|\alpha(x)|+A(x))\leq M.
\eeqnn
For $f\in \mathcal{B}_b(E)$ and $(t,x)\in (0,\infty)\times E$, define
 \beqlb\label{1.4}
 T_tf(x)=\Pi_x\left[e^{\int_0^t\alpha(\xi_s)ds}f(\xi_t)\right].
 \eeqlb
It is well known that $T_tf(x)=\mathbb{E}_{\delta_x}[\langle f,X_t\rangle]$ for every $x\in E$. Then (\ref{1.1})
can be written into
\beqlb\label{1.5}
V_t(f)(x)+\int_0^t\int_E\psi_0(y,V_{t-s}f(y))\beta(y)T_s(x,dy)ds=T_tf(x),
\eeqlb
where $\psi_0(x,\lambda)=\psi(x,\lambda)+a(x)\lambda$, see Theorem 2.23 in \cite{Li11}.

It is shown in \cite{RSZ14} that there exists a family of continuous strictly positive
symmetric functions $\{q_t(x,y), t>0\}$ on $E\times E$ such that $q_t(x,y)\leq e^{Mt}p_t(x,y)$
and for any $f\in \mathcal{B}_b(E)$,
 \beqnn
T_tf(x)=\int_E q_t(x,y)f(y)m(dy).
 \eeqnn
It follows immediately that, for any $p\geq 1$, $\{T_t:t\geq 0\}$ is a strongly continuous semigroup
on $L^p(E,m)$ and
 \beqnn
\|T_tf\|_p^p\leq e^{Mpt}\|f\|_p^p.
 \eeqnn
Define $a_t(x):=q_t(x,x)$. It follows from the assumptions (a) and (b) in the previous subsection that
$a_t$ enjoys the following properties:
\begin{enumerate}
\item[(i)] For any $t>0$, we have
\beqnn
\int_Ea_t(x)m(dx)<\infty.
\eeqnn

\item[(ii)] There exists $t_0>0$ such that for all $t\geq t_0$, $a_t(x)\in L^2(E,m)$.
\end{enumerate}

It follows from (i) above that, for any $t>0$, $T_t$ is a compact operator. The infinitesimal
generator of $\{T_t: t\geq 0\}$ in $L^2(E,m)$ has purely discrete spectrum with eigenvalues
$-\lambda_1>-\lambda_2>-\lambda_3>\cdots$. The first eigenvalue $-\lambda_1$ is simple and
the eigenfunction $\phi_1$ associated with  $-\lambda_1$ can be chosen to be strictly positive
everywhere and continuous. We will assume that $\|\phi_1\|_2=1$. $\phi_1$ is sometimes denoted
as $\phi^{(1)}_1$. For $k>1$, let $\{\phi^{(k)}_j, j=1,2,\cdots n_k\}$ be an orthogonal basis
of the eigenspace (which is finite dimensional) associated with $-\lambda_k$. It is well-known that
$\{\phi^{(k)}_j, j=1,2,\cdots n_k; k=1,2,\ldots\}$ forms a complete orthogonal basis of $L^2(E,m)$
and all the eigenfunctions are continuous. For any $k\geq 1$, $j=1,\ldots, n_k$ and $t>0$, we have
$T_t\phi^{(k)}_j(x)=e^{-\lambda_k t}\phi^{(k)}_j(x)$ and
 \beqnn
 e^{-\lambda_k t/2}|\phi^{(k)}_j(x)|\leq a_t(x)^{1/2},~~x\in E.
 \eeqnn
It follows from the relation above that all the eigenfunctions $\phi^{(k)}_j$ belong to $L^4(E,m)$.
For any $x,y\in E$ and $t>0$, we have
 \beqnn
q_t(x,y)=\sum_{k=1}^\infty e^{-\lambda_k t}\sum_{j=1}^{n_k}\phi^{(k)}_j(x)\phi^{(k)}_j(y),
 \eeqnn
where the series is locally uniformly convergent on $E\times E$. The basic facts recalled in this
paragraph are well known, for instance, one can refer to (\cite{DS84}, Section 2).
In this paper, since we assume that the superprocess $X$ is supercritical, $\lambda_1<0$.

\section{Preliminaries}

\setcounter{equation}{0}

We will use $\<\cdot,\cdot\>_m$ to denote the inner product in $L^2(E,m)$. Any $f\in L^2(E,m)$ admits the following
expansion:
 \beqlb\label{2.1}
f(x)=\sum_{k=1}^\infty\sum_{j=1}^{n_k}a_j^k\phi_j^{(k)}(x),
 \eeqlb
where $a_j^k=\<f,\phi_j^{(k)}\>_m$ and the series converges in $L^2(E,m)$. $a^1_1$ will sometimes be written as $a_1$.
For $f\in L^2(E,m)$, define
\[
\gamma(f):=\inf\{k\geq 1: \mbox{there exists}~j~\mbox{with}~1\leq j\leq n_k~\mbox{such that}~a_j^k\neq 0\},
\]
where we use the usual convention $\inf\emptyset=\infty$. Define
\beqnn
f_1(x)\ar:=\ar\sum_{j=1}^{n_{\gamma(f)}}a_j^{\gamma(f)}\phi_j^{(\gamma(f))}(x).
\eeqnn
We note that if $f\in L^2(E,m)$ is nonnegative and
$m(x: f(x)>0)>0$, then $\<f,\phi_1\>_m >0$ which implies $\gamma(f)=1$. Define
 \beqnn
\mathcal{C}_l\ar:=\ar\left\{g(x)=\sum_{2\lambda_k<\lambda_1}\sum_{j=1}^{n_k}a_j^k\phi_j^{(k)}(x): a_j^k\in \mathbb{R}~\mbox{and}~g\neq 0\right\},\\
\mathcal{C}_c\ar:=\ar\left\{g(x)=\sum_{j=1}^{n_k}a_j^k\phi_j^{(k)}(x): 2\lambda_k=\lambda_1, a_j^k\in \mathbb{R}~\mbox{and}~g\neq 0\right\}
 \eeqnn
and
 \beqnn
\mathcal{C}_s:=\left\{g(x)\in L^2(E,m)\cap L^4(E,m): g\neq 0~\mbox{and}~\lambda_1<2\lambda_{\gamma(g)}\right\}.
 \eeqnn
Note that $\mathcal{C}_l$ consists of these functions in $L^2(E,m)\cap L^4(E,m)$ that only have nontrivial
projection onto the eigen-spaces corresponding to those ``large" eigenvalues $-\lambda_k$ satisfying
$\lambda_1>2\lambda_k$. The space $\mathcal{C}_l$ is of finite dimension. The space $\mathcal{C}_c$ is the (finite dimensional) eigenspace
corresponding to the ``critical" eigenvalue $-\lambda_k$ with $\lambda_1=2\lambda_k$. Note that there may not
be a critical eigenvalue and in this case, $\mathcal{C}_c$ is empty. The space $\mathcal{C}_s$ consists of these functions in
$L^2(E,m)\cap L^4(E,m)$ that only have nontrivial projections onto the eigen-spaces corresponding
to those ''small" eigenvalues $-\lambda_k$ satisfying $\lambda_1<2\lambda_k$. The space $\mathcal{C}_s$ is of infinite dimensional in
general.

We use $\<f,\nu\>:=\int_{\mathbb{R}}f(x)\nu(dx)$.  And whenever we deal with an initial configuration $\mu\in M(E)$, we are implicitly
assuming that it has compact support.

\subsection{Excursion measures of $\{Y_t,t\geq 0\}$}

Let $\mathbb{D}$ denote the space of all paths
$\{w_t:t\geq 0\}$ from $[0,\infty)$ to $M(E)$ that are right continuous in $M(E)$ having zero as a trap. Let $(\mathcal{A}, \mathcal{A}_t)$
denote the natural $\sigma$-algebras on $\mathbb{D}$ generated by the coordinate process.
It is known from \cite[Chapter 8]{Li11} that one can associate with $\{\mathbb{P}_{\delta_x}: x\in E\}$ a family of
$\sigma$-finite measures
$\{\mathbb{N}_x:x\in E\}$, defined on $(\mathbb{D}, \mathcal{A})$ such that $\mathbb{N}_x(\{0\})=0$,
\beqlb\label{2.2}
\mathbb{N}_x(1-e^{-\<f,w_t\>})=-\log\mathbb{P}_{\delta_x}(e^{-\<f,X_t\>}),~f\in \mathcal{B}^+_b(E),~t>0,
\eeqlb
and for every $0<t_1<\cdots<t_n<\infty$, and nonzero $\mu_1, \cdots, \mu_n\in M(E)$,
\beqnn
\mathbb{N}_x(w_{t_1}\in d\mu_1,\cdots,w_{t_n}\in d\mu_n)=\mathbb{N}_x(w_{t_1}\in d\mu_1)\mathbb{P}_{\mu_1}(X_{t_2-t_1}\in d\mu_2)
\cdots\mathbb{P}_{\mu_{n-1}}(X_{t_n-t_{n-1}}\in d\mu_n).
\eeqnn
For earlier work on excursion measures of superprocesses, see \cite{D04,EN91,Li03}.
Next we list some properties of $\mathbb{N}_x$ which will be used later.

\begin{proposition}\label{p2.1}
If $\mathbb{P}_{\delta_x}|\<f,X_t\>|<\infty$, then
\beqnn
\int_\mathbb{D}\<f,w_t\>\mathbb{N}_x(dw)=\mathbb{P}_{\delta_x}\<f,X_t\>.
\eeqnn
If $\mathbb{P}_{\delta_x}\<f,X_t\>^2<\infty$, then
\beqnn
\int_\mathbb{D}\<f,w_t\>^2\mathbb{N}_x(dw)=\mathbb{V}\mbox{ar}_{\delta_x}\<f,X_t\>.
\eeqnn
\end{proposition}

Assume that Condition 8.5 in Li \cite{Li11} holds:
there is a spatially constant local branching mechanism $z\mapsto \psi_\ast(z)$
so that $\psi_\ast'(z)\rightarrow\infty$ as $z\rightarrow\infty$ and $\psi$ is bounded below by $\psi_\ast(z)$ in the sense
\beqnn
\psi(x,f)\geq \psi_\ast(f(x)), ~x\in E, ~f\in \mathcal{B}^+_b(E),
\eeqnn
then it is sufficient for the cumulant semigroup
$(V_t)_{t\geq0}$ to admit the following representation for all $x\in E$:
\beqnn
V_tf(x) = \lambda_t(x, f)+\int_{M(E)^\circ}(1-e^{-\nu(f)})L_t(x,d\nu),
\eeqnn
where $\lambda_t(x,dy)$ is a bounded kernel on $E$ and $(1\wedge \nu(1))L_t(x,d\nu)$ is a bounded
kernel from $E$ to $M(E)^\circ$.
Let
\beqnn
\gamma_t=\int_0^t\int_E\eta(dx)\lambda_s(x,\cdot)ds, H_t=\int_E\eta(dx)L_t(x,\cdot)+HQ_t^\circ~\mbox{and}
~G_t=\int_0^tH_sds, ~t\geq0,
\eeqnn
where $(Q_t^\circ)_{t\geq0}$ denotes the restriction of $(Q_t)_{t\geq0}$ to $M(E)^\circ$.
By Li \cite[Theorems 9.5-9.6]{Li11}, $(G_t)_{t\geq 0}$ is an entrance rule for $(Q_t^\circ)_{t\geq 0}$
and has the decomposition
\beqnn
G_t=\int_0^tG_{t-s}^s\zeta(ds),~~~t\geq 0,
\eeqnn
where $\zeta(ds)$ is a diffuse Radon measure on $[0,\infty)$ and $\{(G_{t}^s)_{t>0}: s\geq 0\}$ is a
family of entrance laws for $(Q_t^\circ)_{t\geq 0}$. Thus
by Li \cite[Theorem A.40]{Li11}, there corresponds a $\sigma$-finite
measure $\mathbf{Q}_s(dw)$ on $(\mathbb{D}, \mathcal{A})$ such that
for every $0<t_1<\cdots<t_n<\infty$, and nonzero $\mu_1, \cdots, \mu_n\in M(E)$,
\beqnn
\mathbf{Q}_s(w_{t_1}\in d\mu_1,\cdots,w_{t_n}\in d\mu_n)=G^s_{t_1}(d\mu_1)Q^\circ_{t_2-t_1}(\mu_1,d\mu_2)
\cdots Q^\circ_{t_n-t_{n-1}}(\mu_{n-1},d\mu_n).
\eeqnn
Let $\mathcal{F}_t=\sigma(X_s: s\in [0, t])$.
Suppose that $N(ds,dw)$ is a Poisson random measure on $(0,\infty)\times \mathbb{D}$
independent of $\{(X_t, \mathcal{F}_t): t\geq 0\}$ with intensity $\zeta(ds)\mathbf{Q}_s(dw)$,
in a probability space $(\widetilde{\Omega},\widetilde{\mathcal{F}},\mathbf{P})$. Define
another process $\{\Lambda_t: t\geq0\}$ by
 \beqnn
\Lambda_t:=\gamma_t+\int_{(0,t]}\int_{\mathbb{D}}w_{t-s}N(ds,dw),~~t\geq 0.
 \eeqnn
The Markov property of $(\Lambda_t)_{t\geq 0}$
follows from a similar proof as that of Li \cite[Theorem 9.29]{Li11}. Moreover,
\beqnn
\mathbf{P}[\exp\{-\Lambda_t(f)\}]\ar=\ar\exp\left\{-\gamma_t(f)-\int_0^t\zeta(ds)\int_{M(E)^\circ}(1-e^{-\<f,\mu\>})G^s_{t-s}(d\mu) \right\}\\
\ar=\ar\exp\left\{-\int_0^t\eta(V_sf)ds-\int_0^t\int_{M(E)^\circ}\mathbb{P}_\mu[1-\exp\{-\<f,X_{t-s}\>\}] H(d\mu)ds\right\}\\
\ar=\ar\exp\left\{-\int_0^t\eta(V_sf)ds-\int_0^t\int_{M(E)^\circ}(1-\exp\{-\<\mu,V_sf\>\}) H(d\mu)ds\right\}\\
\ar=\ar\exp\left\{-\int_0^t\varphi(V_sf)ds\right\}.
\eeqnn
Let $\widetilde{\mathcal{F}}_t$ be the $\sigma$-algebra generated by random variables
$\{N^H(A): A\in \mathcal{B}([0,t])\times\mathcal{A}_t\}$ and
\beqnn
\mathcal{Y}_t=X_t+\Lambda_t,~~X_0=\mu~~\mbox{and}~~\mathcal {H}_t=\sigma(\mathcal{F}_t\cup \widetilde{\mathcal{F}}_t).
\eeqnn
Then $\{(\mathcal{Y}_t,\mathcal{H}_t): t\geq 0\}$
is an superprocess with immigration whose transition semigroup is determined by $(Q_t^N)_{t\geq 0}$.
That is, $\{\mathcal{Y},(\mathcal{H}_t)_{t\geq 0}, \mathbf{P}_\mu\}$ has the same law as
$\{Y,(\mathscr{G}_t)_{t\geq 0}, \mathbb{P}_\mu\}$.

\smallskip

Notice that, for $f\in L^2(E,m)$, $\mathbb{N}_x(\<|f|,\omega_t\>)=T_t|f|(x)<\infty$, which implies that
$\mathbb{N}_x(\<|f|,\omega_t\>=\infty)=0$. Thus, for $f\in L^2(E,m)$,
\beqnn
\mathbb{P}_\mu\left(e^{i\theta\<f,Y_t\>}\right)\ar=\ar\exp\left\{\int_E\int_\mathbb{D}(e^{i\theta \<f,w_t\>}-1)\mathbb{N}_x(dw)\mu(dx)\right\}\\
\ar \ar~~+\exp\left\{i\theta\gamma_t(f)+\int_0^t\zeta(ds)\int_{M(E)^\circ}(\exp\{i\theta \<f,\mu\>\}-1)G^s_{t-s}(d\mu)\right\}.
\eeqnn
Thus, by the Markov property of the process, we have
\beqlb\label{2.3} \nonumber
\ar \ar\mathbb{P}_{\mu}[\exp\{i\theta\<f,Y_{t+s}\>\}|Y_t]\\\nonumber
\ar \ar=\exp\left\{\int_E\int_\mathbb{D}(e^{i\theta \<f,w_s\>}-1)\mathbb{N}_x(dw)Y_t(dx)\right\}\\\nonumber
\ar \ar~~+\exp\left\{i\theta\gamma_s(f)+\int_t^{t+s}\zeta(du)\int_{M(E)^\circ}(\exp\{i\theta \<f,\mu\>\}-1)G^u_{t+s-u}(d\mu)\right\}\\\nonumber
\ar \ar=\exp\left\{\int_E\int_\mathbb{D}(e^{i\theta \<f,w_s\>}-1)\mathbb{N}_x(dw)Y_t(dx)\right\}\\
\ar \ar~~+\exp\left\{\int_0^s\eta(V_u(i\theta f))du+\int_0^s\int_{M(E)^\circ}\mathbb{P}_\mu(\exp\{i\theta \<f,X_{s-u}\>\}-1)H(d\mu)du\right\}.
\eeqlb
The intuitive meaning of (\ref{2.3}) is clear, given $Y_t$,
the population at time $t+s$ is made up of two parts;
the native part generated by the mass $Y_t$ and the
immigration in the time interval $(t, t+s]$.

\subsection{Estimates on the moments of  $\{Y_t,t\geq 0\}$}

\noindent For reader's convenience, we first recall some results about the semigroup $(T_t)$,
the proofs of which can be found in \cite{RSZ14}. For two positive functions $f$ and
$g$ on $E$, by $f(x)\lesssim g(x)$, we will denote the fact that there exists a constant $c>0$ such that $f(x)\leq cg(x)$ for
any $x\in E$ (whose exact value is not relevant to following calculations).

\begin{lemma}\label{l2.1} \cite[Lemma 2.4]{RSZ14}
For any $f\in L^2(E,m)$, $x\in E$ and $t>0$, we have
 \beqlb\label{2.4}
T_tf(x)=\sum_{k=\gamma(f)}^\infty e^{-\lambda_k t}\sum_{j=1}^{n_k}a_j^k\phi_j^{(k)}(x)
 \eeqlb
and
 \beqlb\label{2.5}
\lim_{t\rightarrow\infty}e^{\lambda_{\gamma(f)}t}T_tf(x)=\sum_{j=1}^{n_{\gamma(f)}}a_j^{\gamma(f)}\phi_j^{(\gamma(f))}(x),
 \eeqlb
where the series in (\ref{2.4}) converges absolutely and uniformly in any compact subset of $E$. Moreover,
for any $t_1>0$,
\beqlb\label{2.6}
\ar \ar\sup_{t>t_1}e^{\lambda_{\gamma(f)}t}|T_tf(x)|\leq e^{\lambda_{\gamma(f)}t_1}\|f\|_2\left(\int_Ea_{t_1/2}(x)m(dx)\right)(a_{t_1}(x))^{1/2},\\\nonumber
\ar \ar\sup_{t>t_1}e^{(\lambda_{\gamma(f)+1}-\lambda_{\gamma(f)})t}|e^{\lambda_{\gamma(f)}t}T_tf(x)-f_1(x)|\\ \nonumber
\ar \ar~~\leq e^{\lambda_{\gamma(f)+1}t_1}\|f\|_2\left(\int_Ea_{t_1/2}(x)m(dx)\right)(a_{t_1}(x))^{1/2}.
\eeqlb
\end{lemma}

\begin{lemma}\label{l2.1} \cite[Proposition 9.14]{Li11} Suppose that $\int_{M(E)^\circ}\nu(1)^2H(d\nu)<\infty$.
Then for $t\geq 0$, $\mu\in M(E)$ and
$f\in\mathcal{B}_b(E)$ we have
\beqlb\label{2.7}
\mathbb{P}_\mu\<Y_t,f\>=\mu(T_tf)+\int_0^t\Gamma(T_sf)ds
\eeqlb
and
\beqlb\label{2.8}
\nonumber\mathbb{P}_\mu\<Y_t,f\>^2\ar=\ar\left(\mu(T_tf)+\int_0^t\Gamma(T_sf)ds\right)^2
+\int_E\int_0^t T_s[A(T_{t-s}f)^2](x)ds\mu(dx)\\ \nonumber
\ar \ar+\int_0^t\int_E\int_0^u T_s[A(T_{u-s}f)^2](x)ds\Gamma(dx)du\\
\ar \ar+\int_0^t\int_{M(E)^\circ}\nu(T_sf)^2H(d\nu)ds,
\eeqlb
where $t\rightarrow T_tf$ is defined by (\ref{1.4}) and $\Gamma(f)=\eta(f)+\int_{M(E)^\circ}\nu(f)H(d\nu)$.
\end{lemma}

Let
\beqnn
m(t)\ar=\ar\int_0^{t}\Gamma(a_{2s}^{1/2})ds~\mbox{and}~
n(t)=\int_0^{t}\int_{M(E)^\circ}\left[\nu\left(a_{2s}^{1/2}\right)\right]^2H(d\nu)ds.
\eeqnn
\noindent{\bf Assumption 2.1}  We assume throughout this paper that
there exists $t_0>0$ such that
\beqnn
\int_0^{t_0}\Gamma(a_s^{1/2})ds<\infty~\mbox{and}~\int_0^{t_0}\int_{M(E)^\circ}[\nu(a_s^{1/2})]^2H(d\nu)ds<\infty.
\eeqnn
\noindent{\bf Remark 2.1} It follows from Assumption 2.1 and the fact $t\rightarrow a_t$ is decreasing that, for any $t>0$,
\beqnn
\Gamma(a_t^{1/2})<\infty.
\eeqnn
Furthermore, for any $t\geq t_0$, $m(t)<\infty$, $n(t)<\infty$.
Since $e^{-\lambda_k t/2}|\phi^{(k)}_j(x)|\leq a_t(x)^{1/2}$, we have
\beqnn
|\Gamma(\phi^{(k)}_j)|\leq e^{\lambda_k t/2}|\Gamma(a_t^{1/2})|<\infty
\eeqnn
for any $k\geq 1$, $j=1,\ldots, n_k$.

\medskip

For $f\in L^2(E,m)\cap L^4(E,m)$, and $x\in E$, it follows from H\"{o}lder's inequality that
\beqlb\label{2.9}
(T_{t-s}f)^2(x)\leq e^{M(t-s)}T_{t-s}(f^2)(x).
\eeqlb
Thus, using a routine limit argument, one can easily check that (\ref{2.7}) and (\ref{2.8}) also
hold for $f\in L^2(E,m)\cap L^4(E,m)$ under Assumption 2.1.

\begin{lemma}
Under Assumption 2.1, for any $f\in L^2(E,m)\cap L^4(E,m)$,
\begin{enumerate}
\item[(1)]If $\lambda_1< 2\lambda_{\gamma(f)}$, then for any $x\in E$,
 \beqlb\label{2.10}
\lim_{t\rightarrow\infty}e^{\lambda_1t/2}\mathbb{P}_{\delta_x}\<f,Y_t\>=0
 \eeqlb
and
 \beqlb\label{2.11}\nonumber
\lim_{t\rightarrow\infty}e^{\lambda_1t}\mathbb{P}_{\delta_x}\<f,Y_t\>^2\ar=\ar
\int_0^\infty e^{\lambda_1s}
\<A(T_sf)^2,\phi_1\>_m ds\phi_1(x)\\
\ar \ar+\int_0^\infty\int_l^\infty e^{\lambda_1l}
\<A(T_{u-l}f)^2,\phi_1\>_m dudl \Gamma(\phi_1).
 \eeqlb

\item[(2)]If $\lambda_1= 2\lambda_{\gamma(f)}$, then for any $(t,x)\in (3t_0,\infty)\times E$,
\beqlb\label{2.12}
\lim_{t\rightarrow\infty}t^{-1}e^{\lambda_1t}\mathbb{V}\mbox{ar}_{\delta_x}\<f,Y_t\>=
\left(\phi_1(x)+\frac{\Gamma(\phi_1)}{-\lambda_1}\right)\<Af_1^2,\phi_1\>_m.
\eeqlb

\item[(3)]If $\lambda_1> 2\lambda_{\gamma(f)}$, then for any $x\in E$,
 \beqlb\label{2.13}\nonumber
\lim_{t\rightarrow\infty}e^{2\lambda_{\gamma(f)}t}\mathbb{V}\mbox{ar}_{\delta_x}\<f,Y_t\>
\ar=\ar\int_0^\infty e^{2\lambda_{\gamma(f)}s}T_s(Af_1^2)(x)ds\\\nonumber
\ar \ar+\int_0^\infty\int_E
\int_l^\infty e^{2\lambda_{\gamma(f)}l}T_{u-l}[A(f_1)^2](x)du\Gamma(dx)dl\\
\ar \ar+\int_0^\infty\int_{M(E)^\circ} e^{2\lambda_{\gamma(f)}s}\nu(f_1^2)(x)H(d\nu)ds.
\eeqlb
\end{enumerate}
\end{lemma}

\noindent{\it Proof.} (1) If $\lambda_1< 2\lambda_{\gamma(f)}$, by the moment formula (\ref{2.7}),
\beqnn
e^{\lambda_1t/2}|\mathbb{P}_{\delta_x}\<f,Y_t\>|\ar\leq\ar e^{(\lambda_1- 2\lambda_{\gamma(f)})t/2}(e^{\lambda_{\gamma(f)}t}|T_tf(x)|)
+e^{\lambda_1t/2}\int_0^t|\Gamma(T_sf)|ds.
\eeqnn
Using H\"{o}lder's inequality, we can get,
\beqlb\label{2.14}
|T_sf(x)|\leq \int_E q_s(x,y)|f|(y)m(dy)\leq
\|f\|_2 \left(\int_E q_s(x,y)^2m(dy)\right)^{1/2}= \|f\|_2a_{2s}(x)^{1/2}.
\eeqlb
Then for any $t>t_0>0$, we have by (\ref{2.6}),
\beqlb\label{2.15}\nonumber
e^{\lambda_1t/2}\int_0^t\Gamma(T_sf)ds\ar=\ar e^{\lambda_1t/2}\int_0^{t_0}\Gamma(T_sf)ds+e^{\lambda_1t/2}\int_{t_0}^t\Gamma(T_sf)ds\\\nonumber
\ar\leq\ar e^{\lambda_1t/2}\|f\|_2\int_0^{t_0}\Gamma(a_{2s}^{1/2})ds +\int_{t_0}^te^{\lambda_1(t-s)/2}e^{(\lambda_1- 2\lambda_{\gamma(f)})s/2}\Gamma(e^{\lambda_{\gamma(f)}s}T_sf)ds\\
\ar\lesssim\ar e^{\lambda_1t/2}\|f\|_2m(t_0)+[e^{(\lambda_1- 2\lambda_{\gamma(f)})t/2}+e^{\lambda_1t/2}]\Gamma(a_{t_0}^{1/2}).
\eeqlb
Combining the estimates above, we get
\beqnn
\lim_{t\rightarrow\infty}e^{\lambda_1t/2}\mathbb{P}_{\delta_x}\<f,Y_t\>=0.
 \eeqnn
For the second moment, we have already known from the proof of Lemma 2.3 in \cite{RSZ14} that
 \beqnn
\lim_{t\rightarrow\infty}e^{\lambda_1t}\int_0^t T_s[A(T_{t-s}f)^2](x)ds=\int_0^\infty e^{\lambda_1s}
\<A(T_sf)^2,\phi_1\>_m ds\phi_1(x)
 \eeqnn
and for $s>3t_0$,
 \beqlb\label{2.16}
e^{\lambda_1s}\int_0^s T_u[A(T_{s-u}f)^2](x)du\lesssim a_{t_0}(x)^{1/2}.
 \eeqlb
For any $s\leq 3t_0$, using (\ref{2.9}) and (\ref{2.14}), we get
\beqnn
\int_E \int_0^s T_u[A(T_{s-u}f)^2](x)du\Gamma(dx)
\ar\leq\ar\int_E \int_0^s Me^{M(s-u)}T_s(f^2)(x)du\Gamma(dx)\\
\ar\leq\ar e^{3Mt_0} \|f^2\|_2\Gamma(a_{2s}^{1/2}).
 \eeqnn
Combining with (\ref{2.16}), we get
\beqlb\label{2.17}      \nonumber
\ar \ar e^{\lambda_1t}\int_0^t\int_E\int_0^s T_u[A(T_{s-u}f)^2](x)du\Gamma(dx)ds\\ \nonumber
\ar=\ar\int_{3t_0}^te^{\lambda_1(t-s)}\int_Ee^{\lambda_1s}
\int_0^s T_u[A(T_{s-u}f)^2](x)du\Gamma(dx)ds\\ \nonumber
\ar \ar+ e^{\lambda_1t}\int_0^{3t_0}\int_E\int_0^s T_u[A(T_{s-u}f)^2](x)du\Gamma(dx)ds\\
\ar\lesssim\ar\Gamma(a_{t_0}^{1/2})+e^{\lambda_1t}e^{3Mt_0}\|f^2\|_2m(3t_0).
 \eeqlb
Thus, by the dominated convergence theorem and let $l=t-s$, we have,
 \beqnn
\ar \ar e^{\lambda_1t}\int_0^t\int_E\int_0^s T_{s-u}[A(T_uf)^2](x)du\Gamma(dx)ds\\
\ar=\ar \int_0^te^{\lambda_1(t-s)}\int_Ee^{\lambda_1s}\int_0^s T_{s-u}[A(T_uf)^2](x)du\Gamma(dx)ds\\
\ar=\ar \int_0^te^{\lambda_1l}\int_Ee^{\lambda_1(t-l)}\int_0^{t-l}T_{t-l-u}[A(T_uf)^2](x)du\Gamma(dx)dl\\
\ar=\ar \int_0^t\int_Ee^{\lambda_1l}e^{\lambda_1(t-l)}\int_l^tT_{t-u}[A(T_{u-l}f)^2](x)du\Gamma(dx)dl\\
\ar\rightarrow\ar\int_0^\infty\int_l^\infty e^{\lambda_1l}
\<A(T_{u-l}f)^2,\phi_1\>_m\Gamma(\phi_1) dudl,~~~\mbox{as}~t\rightarrow\infty.
 \eeqnn
For the last term in (\ref{2.8}), for $t> t_0$, we have by (\ref{2.14}),
 \beqlb\label{2.18} \nonumber
\ar \ar e^{\lambda_1t}\int_0^t\int_{M(E)^\circ}\nu(T_sf)^2H(d\nu)ds\\ \nonumber
\ar \ar=\int_{t_0}^t e^{\lambda_1(t-s)}\int_{M(E)^\circ}e^{(\lambda_1- 2\lambda_{\gamma(f)})s}\nu(e^{\lambda_{\gamma(f)}s}T_sf)^2H(d\nu)ds
+e^{\lambda_1t}\int_0^{t_0}\int_{M(E)^\circ}\nu(T_sf)^2H(d\nu)ds\\
\ar \ar\lesssim (e^{(\lambda_1- 2\lambda_{\gamma(f)})t}+e^{\lambda_1t})\int_{M(E)^\circ}\nu\left(a_{t_0}^{1/2}\right)^2H(d\nu)
+ e^{\lambda_1t}\|f\|^2_2 n(t_0),
 \eeqlb
so we have
\beqnn
\lim_{t\rightarrow\infty}e^{\lambda_1t}\int_0^t\int_{M(E)^\circ}\nu(T_sf)^2H(d\nu)ds=0.
\eeqnn
(2) If $\lambda_1= 2\lambda_{\gamma(f)}$,
\beqnn
t^{-1}e^{\lambda_1t}\mathbb{V}\mbox{ar}_{\delta_x}\<f,Y_t\>\ar=\ar t^{-1}e^{\lambda_1t}\int_0^t T_s[A(T_{t-s}f)^2](x)ds\\
\ar \ar+~t^{-1}e^{\lambda_1t}\int_0^t\int_E\int_0^u T_s[A(T_{u-s}f)^2](x)ds\Gamma(dx)du\\
\ar \ar+~t^{-1}e^{\lambda_1t}\int_0^t\int_{M(E)^\circ}\nu(T_sf)^2H(d\nu)ds\\
\ar:=\ar A_1(t,x)+A_2(t)+A_3(t).
\eeqnn
We have already known from the proof of Lemma 2.3 in \cite{RSZ14} that for $t>3t_0$,
\beqnn
|A_1(t,x)-\<Af_1^2,\phi_1\>_m\phi_1(x)|\lesssim t^{-1}a_{t_0}(x)^{1/2}.
\eeqnn
Thus, by similar estimate to (\ref{2.17}), we have
\beqnn
A_2(t)\ar=\ar\int_0^tt^{-1}e^{\lambda_1(t-s)}\int_Ee^{\lambda_1s}\int_0^s T_u[A(T_{s-u}f)^2](x)du\Gamma(dx)ds\\
\ar=\ar\int_{3t_0}^t\frac{s}{t}e^{\lambda_1(t-s)}s^{-1}e^{\lambda_1s}\int_E\int_0^s T_u[A(T_{s-u}f)^2](x)du\Gamma(dx)ds\\
\ar \ar +~t^{-1}e^{\lambda_1t}\int_0^{3t_0}\int_E\int_0^s T_u[A(T_{s-u}f)^2](x)du\Gamma(dx)ds\\
\ar\lesssim\ar t^{-1}\Gamma(a_{t_0}^{1/2})+\<Af_1^2,\phi_1\>_m\Gamma(\phi_1)+t^{-1}e^{\lambda_1t}\|f^2\|_2m(3t_0).
\eeqnn
Let $l=t-s$, by the dominated convergence theorem, we have as $t\rightarrow\infty$,
\beqnn
A_2(t)\ar=\ar\int_0^te^{\lambda_1l}\frac{t-l}{t}\int_E(t-l)^{-1}e^{\lambda_1(t-l)}\int_l^tT_{t-u}[A(T_{u-l}f)^2](x)du\Gamma(dx)dl \\
\ar\rightarrow\ar \int_0^\infty e^{\lambda_1l}\Gamma(\phi_1)\<Af_1^2,\phi_1\>_mdl
\eeqnn
and similarly to (\ref{2.18}), we can prove that
\beqnn
A_3(t)\rightarrow 0~\mbox{as}~t\rightarrow\infty.
\eeqnn

\noindent(3) If $\lambda_1> 2\lambda_{\gamma(f)}$, we have already known from the proof of Lemma 2.3 in \cite{RSZ14} that
\beqnn
\lim_{t\rightarrow\infty}e^{2\lambda_{\gamma(f)}t}\int_0^tT_s[A(T_{t-s}f)^2](x)ds=\int_0^\infty e^{2\lambda_{\gamma(f)}s}T_s(Af_1^2)(x)ds.
\eeqnn
For the second term in (\ref{2.8}), by (2.37) in \cite{RSZ14}, for $t>3t_0$,
\beqlb\label{2.19} \nonumber
\ar \ar e^{2\lambda_{\gamma(f)}t}\int_0^t\int_E\int_0^s T_u[A(T_{s-u}f)^2](x)du\Gamma(dx)ds\\ \nonumber
\ar=\ar\int_{3t_0}^te^{2\lambda_{\gamma(f)}(t-s)}\int_Ee^{2\lambda_{\gamma(f)}s}
\int_0^s T_u[A(T_{s-u}f)^2](x)du\Gamma(dx)ds\\ \nonumber
\ar \ar+\int_0^{3t_0}e^{2\lambda_{\gamma(f)}t}\int_E\int_0^s T_u[A(T_{s-u}f)^2](x)du\Gamma(dx)ds\\
\ar\lesssim\ar\Gamma(a_{t_0}^{1/2})+e^{2\lambda_{\gamma(f)}t}e^{3Mt_0}\|f^2\|_2m(3t_0).
 \eeqlb
Consequently, we can use the dominated convergence theorem to get
\beqnn
\ar \ar e^{2\lambda_{\gamma(f)}t}\int_0^t\int_E\int_0^s T_u[A(T_{s-u}f)^2](x)du\Gamma(dx)ds\\
\ar=\ar\int_0^t\int_E
e^{2\lambda_{\gamma(f)}l}e^{2\lambda_{\gamma(f)(t-l)}}\int_l^t T_{u-l}[A(T_{t-u}f)^2](x)du\Gamma(dx)dl\\
\ar\rightarrow\ar\int_0^\infty\int_E
\int_l^\infty e^{2\lambda_{\gamma(f)}l}T_{u-l}[A(f_1)^2](x)du\Gamma(dx)dl.
 \eeqnn
For the last term in (\ref{2.8}),
\beqlb\label{2.20}\nonumber
e^{2\lambda_{\gamma(f)}t}\int_0^t\int_{M(E)^\circ}\nu(T_sf)^2H(d\nu)ds
\ar=\ar\int_{t_0}^t e^{2\lambda_{\gamma(f)}(t-s)}\int_{M(E)^\circ}\nu(e^{\lambda_{\gamma(f)}s}T_sf)^2H(d\nu)ds\\\nonumber
 \nonumber\ar \ar+~e^{2\lambda_{\gamma(f)}t}\int_0^{t_0}\int_{M(E)^\circ}\nu(T_sf)^2H(d\nu)ds\\
\ar\lesssim \ar \int_{M(E)^\circ}\nu\left(a_{t_0}^{1/2}\right)^2H(d\nu)+e^{2\lambda_{\gamma(f)}t}\|f\|^2_2n(t_0).
 \eeqlb
Then it follows by the dominated convergence theorem that as $t\rightarrow\infty$,
\beqnn
e^{2\lambda_{\gamma(f)}t}\int_0^t\int_{M(E)^\circ}\nu(T_sf)^2H(d\nu)ds
\ar=\ar\int_0^t e^{2\lambda_{\gamma(f)}(t-s)}\int_{M(E)^\circ}\nu(e^{\lambda_{\gamma(f)}s}T_sf)^2H(d\nu)ds\\
\ar\rightarrow\ar\int_0^\infty\int_{M(E)^\circ} e^{2\lambda_{\gamma(f)}s}\nu(f_1)^2H(d\nu)ds.
\eeqnn
Combining the above, we get
\beqnn
\lim_{t\rightarrow\infty}e^{2\lambda_{\gamma(f)}t}\mathbb{V}\mbox{ar}_{\delta_x}\<f,Y_t\>
\ar=\ar\int_0^\infty e^{2\lambda_{\gamma(f)}s}T_s(Af_1^2)(x)ds\\
\ar \ar+\int_0^\infty\int_E
\int_l^\infty e^{2\lambda_{\gamma(f)}l}T_{u-l}[A(f_1)^2](x)du\Gamma(dx)dl\\
\ar \ar+\int_0^\infty\int_{M(E)^\circ} e^{2\lambda_{\gamma(f)}s}\nu(f_1)^2H(d\nu)ds.
\eeqnn
Thus the proof of (3) is now complete. \qed

\section{Main results}

\setcounter{equation}{0}

\subsection{Law of large numbers of  $\{Y_t,t\geq 0\}$}

Define
\beqnn
W_t^{k,j}\ar:=\ar e^{\lambda_kt}\<\phi_j^{(k)},Y_t\>,\\
H_t^{k,j}\ar:=\ar e^{\lambda_kt}\<\phi_j^{(k)},Y_t\>-
\lambda_k^{-1}(e^{\lambda_kt}-1)\Gamma(\phi_j^{(k)}),~~t\geq 0.
\eeqnn
\begin{lemma}\label{l3.1}
$\{H_t^{k,j}: t\geq0\}$ is a martingale under $\mathbb{P}_\mu$. Moreover, if $\lambda_1>2\lambda_k$,
then $\sup_{t>3t_0}\mathbb{P}_\mu(H_t^{k,j})^2<\infty$. Thus the limit
\beqnn
H_\infty^{k,j}:=\lim_{t\rightarrow\infty}H_t^{k,j}
\eeqnn
exists $\mathbb{P}_\mu$-a.s. and in $L^2(\mathbb{P}_\mu)$.
\end{lemma}

\noindent{\it Proof.} By the moment formula and the fact that $T_t\phi^{(k)}_j(x)=e^{-\lambda_k t}\phi^{(k)}_j(x)$,
\beqnn
\mathbb{P}_\mu[H_{t+s}^{k,j}|\mathscr{G}_t]
\ar=\ar e^{\lambda_kt}\<\phi_j^{(k)},Y_t\>+
e^{\lambda_k(t+s)}\int_t^{t+s}\Gamma(T_{t+s-u}\phi_j^{(k)})du-\lambda_k^{-1}(e^{\lambda_k(t+s)}-1)\Gamma(\phi_j^{(k)})\\
\ar=\ar e^{\lambda_kt}\<\phi_j^{(k)},Y_t\>-
\lambda_k^{-1}(e^{\lambda_kt}-1)\Gamma(\phi_j^{(k)})\\
\ar=\ar H_{t}^{k,j}.
\eeqnn
By (\ref{2.19}) and (\ref{2.20}), we have that
\beqnn
\ar \ar\sup_{t>3t_0}\mathbb{P}_\mu[H_{t}^{k,j}]^2\lesssim \sup_{t>3t_0}e^{2\lambda_k t}\mathbb{V}\mbox{ar}_{\mu}\<\phi_j^{(k)},Y_t\>+\mu(\phi_j^{(k)})^2\\
\ar\lesssim\ar \<a_{t_0}^{1/2},\mu\>+\Gamma(a_{t_0}^{1/2})+\int_{M(E)^\circ}\nu\left(a_{t_0}^{1/2}\right)^2H(d\nu)+m(3t_0)+n(t_0)+\mu(\phi_j^{(k)})^2
\eeqnn
from which the convergence asserted in the lemma follows easily.
\qed

\medskip

\noindent{\bf Remark 3.1} We know from Remark 2.1 that $\Gamma(\phi_j^{(k)})<\infty$. If $\lambda_1>2\lambda_k$,
then
\beqnn
W_\infty^{k,j}:=\lim_{t\rightarrow\infty}e^{\lambda_kt}\<\phi_j^{(k)},Y_t\>=H_\infty^{k,j}-
\lambda_k^{-1}(e^{\lambda_kt}-1)\Gamma(\phi_j^{(k)})
=H_\infty^{k,j}+\lambda_k^{-1}\Gamma(\phi_j^{(k)})\eeqnn
exists $\mathbb{P}_\mu$-a.s. and in $L^2(\mathbb{P}_\mu)$ for any $k\geq 1$, $j=1,\ldots, n_k$.

\begin{theorem}
If $f\in L^2(E,m)\cap L^4(E,m)$ with $\lambda_1>2\lambda_{\gamma(f)}$ and $\Gamma(\phi_j^{(\gamma(f))})<\infty$,
then as $t\rightarrow\infty$,
\beqnn
e^{\lambda_{\gamma(f)}}\<f,Y_t\>\rightarrow \sum_{j=1}^{n_{\gamma(f)}}a_j^{\gamma(f)}W_\infty^{\gamma(f),j},
~~\mbox{in}~L^2(\mathbb{P}_\mu).
\eeqnn
\end{theorem}

\noindent{\it Proof.} By using the moment estimates of $Y$,
the proof is similar to that of \cite[Theorem 1.6]{RSZ14}, we omit the details here. \qed

\medskip

\noindent{\bf Remark 3.2} When $\gamma(f)=1$, let $\widetilde{W}_\infty$ be short for $W_\infty^{1,1}$.
Therefore as $t\rightarrow\infty$,
\beqnn
e^{\lambda_1t}\<f,Y_t\>\rightarrow \<f,\phi_1\>_m\widetilde{W}_\infty,~~\mbox{in}~L^2(\mathbb{P}_\mu).
\eeqnn
In particular, the convergence also holds in $\mathbb{P}_\mu$-probability.

\subsection{Central limit theorems for $\{Y_t,t\geq 0\}$}

\medskip

For $f\in \mathcal{C}_s$ and $h\in \mathcal{C}_c$, we define
\beqnn
\sigma_f^2:=\int_0^\infty e^{\lambda_1s}\<A(T_sf)^2,\phi_1\>_mds
\eeqnn
and
\beqnn
\rho_h^2:=\<Ah^2,\phi_1\>_m.
\eeqnn
For $g(x)=\sum_{k:2\lambda_k<\lambda_1}\sum_{j=1}^{n_k}a_j^k\phi_j^{(k)}(x)\in \mathcal{C}_l$, we define
\beqnn
I_sg(x)=\sum_{k:2\lambda_k<\lambda_1}\sum_{j=1}^{n_k}e^{\lambda_ks}a_j^k\phi_j^{(k)}(x)~~\mbox{and}~~
\beta_g^2:=\int_0^\infty e^{-\lambda_1s}\left\<A(I_sg)^2,\phi_1\right\>_mds.
\eeqnn

\begin{theorem}\label{t3.3} If $f\in \mathcal{C}_s$, $h\in \mathcal{C}_c$ and $g(x)=\sum_{k:2\lambda_k<\lambda_1}
\sum_{j=1}^{n_k}a_j^k\phi_j^{(k)}(x)\in \mathcal{C}_l$,
then $\sigma_f^2<\infty$, $\rho_h^2<\infty$ and $\beta_g^2<\infty$. Furthermore, it holds that,
as $t\rightarrow\infty$,
\beqnn
\ar \ar\left(e^{\lambda_1t}\<\phi_1,Y_t\>,
\frac{\<g,Y_t\>-\sum_{k:2\lambda_k<\lambda_1}\sum_{j=1}^{n_k}e^{-\lambda_kt}a_j^kW_\infty^{k,j}}{\sqrt{\<\phi_1,Y_t\>}},
\frac{\<h,Y_t\>}{\sqrt{t\<\phi_1,Y_t\>}},
\frac{\<f,Y_t\>}{\sqrt{\<\phi_1,Y_t\>}}\right)\\
\ar \ar \stackrel{d}{\rightarrow}(\widetilde{W}_\infty, G_3(g), G_2(h), G_1(f)),
\eeqnn
where
$G_3(g)\sim \mathcal{N}(0,\beta_g^2)$, $G_2(h)\sim \mathcal{N}(0,\rho_h^2)$ and $G_1(f)\sim \mathcal{N}(0,\sigma_f^2)$.
Moreover, $\widetilde{W}_\infty$, $G_3(g)$, $G_2(h)$ and $G_1(f)$ are independent.
\end{theorem}

\noindent{\bf Remark 3.3} In \cite{RSZ15}, Corollaries 1.5, 1.6 and 1.7 are excellent complementary
of Theorem 1.4. They all hold in our situation only with $X_t$ replaced by $Y_t$. We will not restate them here.

\subsection{Proof of the central limit theorem of $Y$}

The general methodology is similar
to that of \cite{RSZ14} and \cite{RSZ15}. We recall some facts about weak convergence which will be
used later. For $f: \mathbb{R}^n\rightarrow \mathbb{R}$, let
$\|f\|_L:=\sup_{x\neq y}|f(x)-f(y)|/\|x-y\|$ and $\|f\|_{BL}:=\|f\|_\infty+\|f\|_L$. For any distributions
$\nu_1$ and $\nu_2$ on $\mathbb{R}^n$, define
\beqnn
d(\nu_1,\nu_2):=\sup\left\{\left|\int fd\nu_1-\int fd\nu_2\right|: \|f\|_{BL}\leq 1\right\}.
\eeqnn
Then $d$ is a metric. It follows from [9, Theorem 11.3.3] that the topology generated by $d$
is equivalent to the weak convergence topology. From the definition, we can easily see that,
if $\nu_1$ and $\nu_2$ are the distributions of two $\mathbb{R}^n$-valued random variables $X$
and $Y$ respectively, then
\beqnn
d(\nu_1,\nu_2)\leq E\|X-Y\|\leq \sqrt{E\|X-Y\|^2}.
\eeqnn
Before the proof of Theorem \ref{t3.3}, we prove several lemmas first. The first lemma below was
proved in \cite{RSZ15}, we state it here for reader's convenience.
Recall the excursion measure $\mathbb{N}_x$ defined by (\ref{2.2}) on the
probability space $(\mathbb{D}, \mathcal{A})$, define
\beqnn
\widetilde{H}_t^{k,j}(w):=e^{\lambda_kt}\<\phi_j^{(k)},w_t\>,~~t\geq 0,~w\in\mathbb{D}.
\eeqnn

\begin{lemma}\label{l3.2} (\cite[Lemma 3.1]{RSZ15}) If $\lambda_1>2\lambda_k$, then the limit
\beqnn
\widetilde{H}_\infty^{k,j}:=\lim_{t\rightarrow\infty}\widetilde{H}_t^{k,j}
\eeqnn
exists $\mathbb{N}_x$-a.e., in $L^1(\mathbb{N}_x)$ and in $L^2(\mathbb{N}_x)$.
\end{lemma}

\begin{lemma}\label{l3.3}
If $f\in \mathcal{C}_s$, then $\sigma^2_f<\infty$ and for any nonzero $\mu\in M(E)$, it holds under
$\mathbb{P}_\mu$ that
\beqnn
\left(e^{\lambda_1t}\<\phi_1, Y_t\>, e^{\lambda_1t/2}\<f,Y_t\> \right)\rightarrow
\left(\widetilde{W}_\infty, G_1(f)\sqrt{\widetilde{W}_\infty}\right),~~t\rightarrow\infty,
\eeqnn
where $G_1(f)\sim \mathcal{N}(0,\sigma^2_f)$. Moreover, $\widetilde{W}_\infty$ and $G_1(f)$ are independent.
\end{lemma}

\noindent{\it Proof.}
We need to consider the limit of the $\mathbb{R}^2$-valued random variable defined by
\beqnn
U_1(t):=\left(e^{\lambda_1t}\<\phi_1,Y_t\>, e^{\lambda_1t/2}\<f,Y_t\>\right),
\eeqnn
or equivalently, we need to consider the limit of $U_1(s+t)$ as $t\rightarrow\infty$ for any $s>0$.
For $s,t>t_0$,
\beqnn
U_1(s+t)\ar=\ar\left(e^{\lambda_1(t+s)}\<\phi_1,Y_{t+s}\>, e^{\lambda_1(t+s)/2}\<f,Y_{t+s}\>-e^{\lambda_1(t+s)/2}\<T_sf,Y_t\>\right)\\
\ar \ar +\left(0,e^{\lambda_1(t+s)/2}\<T_sf,Y_t\>\right).
\eeqnn
We will prove that the second term on the right hand has no contribution to the double limit, first as $t\rightarrow\infty$
and then $s\rightarrow\infty$. The double limit of the first term is equal to another $\mathbb{R}^2$-valued
random variable $U_2(s,t)$ where
\beqnn
U_2(s,t):=\left(e^{\lambda_1t}\<\phi_1,Y_t\>, e^{\lambda_1(t+s)/2}\<f,Y_{t+s}\>-e^{\lambda_1(t+s)/2}\<T_sf,Y_t\>\right).
\eeqnn
We claim that
\beqnn
U_2(s,t)\stackrel{d}{\rightarrow}(\widetilde{W}_\infty, \sqrt{\widetilde{W}_\infty}G_1(s)),
~~\mbox{as}~t\rightarrow\infty,
\eeqnn
where $G_1(s)\sim \mathcal{N}(0, \sigma_f^2(s))$ with $\sigma_f^2(s)$ to be given later. Denote the characteristic
function of $U_2(s,t)$ under $\mathbb{P}_\mu$ by $\kappa(\theta_1,\theta_2,s,t)$:
\beqlb\label{3.1} \nonumber
\kappa(\theta_1,\theta_2,s,t)\ar=\ar\mathbb{P}_\mu\left(\exp\left\{i\theta_1e^{\lambda_1t}\<\phi_1,Y_t\>+
i\theta_2e^{\lambda_1(t+s)/2}\<f,Y_{t+s}\>-i\theta_2e^{\lambda_1(t+s)/2}\<T_sf,Y_{t}\>\right\}\right)\\ \nonumber
\ar=\ar\mathbb{P}_\mu\left(\exp\left\{i\theta_1e^{\lambda_1t}\<\phi_1,Y_t\>\right.\right.\\\nonumber
\ar \ar\left.\left.+\int_E\int_\mathbb{D}\left(\exp\{i\theta_2e^{\lambda_1(t+s)/2}\<f,w_s\>\}-1-i\theta_2e^{\lambda_1(t+s)/2}\<f,w_s\> \right)\mathbb{N}_x(dw)Y_t(dx)\right.\right.\\\nonumber
\ar \ar\left.\left.+\int_0^s\int_{M(E)^\circ}
\mathbb{P}_\mu\left(\exp\{i\theta_2e^{\lambda_1(t+s)/2}\<f,X_{s-u}\>\}-1\right)H(d\mu)du\right.\right.\\
\ar \ar\left.\left.+ \int_0^s\eta(V_u(i\theta_2e^{\lambda_1(t+s)/2} f))du\right\}\right).
\eeqlb
Let
\beqnn
\ar \ar\mathbb{P}_\mu\left(\exp\left\{\int_E\int_\mathbb{D}\left(\exp\{i\theta_2e^{\lambda_1(t+s)/2}\<f,w_s\>\}-1-i\theta_2e^{\lambda_1(t+s)/2}\<f,w_s\> \right)\mathbb{N}_x(dw)Y_t(dx)\right\}\right)\\
\ar \ar=\mathbb{P}_\mu\left(\exp\left\{-\frac{1}{2}\theta^2_2e^{\lambda_1t}\<V_s,Y_t\>
+\<R_s(e^{\lambda_1(t+s)/2}\theta_2,\cdot),Y_t\>\right\}\right),
\eeqnn
where $V_s(x)=e^{\lambda_1s}\mathbb{V}ar_{\delta_x}\<f,X_s\>$ and
\beqnn
R_s(\theta,x)=\int_\mathbb{D}\left(\exp{\<i\theta f,w_s\>}-1- i\theta\<f,w_s\>+\frac{1}{2}\theta^2\<f,w_s\>^2\right)\mathbb{N}_x(dw).
\eeqnn
By Remark 3.2 and the fact that $V_s(x)\lesssim a_{t_0}(x)^{1/2}\in L^2(E,m)\cap L^4(E,m)$,
we have
\beqnn
\lim_{t\rightarrow\infty}e^{\lambda_1t}\<V_s,Y_t\>=\<V_s,\phi_1\>_m \widetilde{W}_\infty,~~\mbox{in probability}.
\eeqnn
Let $Z_s:=e^{\lambda_1s/2}\<f,w_s\>$ and
\beqnn
h(x,s,t):=\mathbb{N}_x\left(Z_s^2\left(\frac{\theta_2e^{\lambda_1t/2}Z_s}{6}\wedge 1\right)\right).
\eeqnn
We know from (2.16) in \cite{RSZ15} that $h(x,s,t)\downarrow 0$ as $t\uparrow \infty$ and for $t>3t_0$,
\beqnn
h(x,s,t)\leq\mathbb{N}_x(Z_s^2)=e^{\lambda_1s}\mathbb{V}\mbox{ar}_{\delta_x}(\<f,X_s\>)\lesssim a_{t_0}(x)^{1/2}\in L^2(E,m).
\eeqnn
By (3.8) in \cite{RSZ15},
\beqnn
|R_s(e^{\lambda_1(t+s)/2}\theta_2,x)|\leq \theta_2^2e^{\lambda_1t}h(x,s,t)\leq \theta_2^2e^{\lambda_1t}a_{t_0}(x)^{1/2}\in L^2(E,m).
\eeqnn
We claim that, as $t\rightarrow\infty$,
\beqlb\label{3.2}
\mathbb{P}_\mu|\<R_s(e^{\lambda_1(t+s)/2}\theta_2,\cdot),Y_t\>|\leq\theta_2^2e^{\lambda_1t}\left[T_t(h(\cdot,s,t))+\int_0^t\Gamma(T_vh)dv\right]
\rightarrow 0.
\eeqlb
For the first term on the right hand side of (\ref{3.2}). By the results in \cite{RSZ15}, for any $u<t$,
\beqnn
\limsup_{t\rightarrow\infty}e^{\lambda_1t}[T_t(h(\cdot,s,t))]\leq \limsup_{t\rightarrow\infty}e^{\lambda_1t}[T_t(h(\cdot,s,u))]
=\<h(\cdot,s,u),\phi_1\>_m\phi_1(x).
\eeqnn
Letting $u\rightarrow\infty$, we get $\lim_{t\rightarrow\infty}e^{\lambda_1t}[T_t(h(\cdot,s,t))]=0$.
For the second term, we have
\beqnn
\limsup_{t\rightarrow\infty}e^{\lambda_1t}\int_0^t\Gamma(T_v(h(\cdot,s,t)))dv
\ar\leq\ar\limsup_{t\rightarrow\infty}\int_0^te^{\lambda_1(t-v)}\Gamma(e^{\lambda_1v}T_v(h(\cdot,s,u)))dv\\
\ar\leq\ar\int_0^\infty e^{\lambda_1v}\<h(\cdot,s,u),\phi_1\>_mdv\Gamma(\phi_1)\\
\ar\lesssim\ar\Gamma(\phi_1)\<(a_{t_0})^{1/2},\phi_1\>_m.
\eeqnn
By the dominated convergence theorem,
\beqnn
\lim_{t\rightarrow\infty}e^{\lambda_1t}\int_0^t\Gamma(T_v(h(\cdot,s,t)))dv
\leq \lim_{u\rightarrow\infty}\lim_{t\rightarrow\infty}\int_0^te^{\lambda_1(t-v)}\Gamma(e^{\lambda_1v}T_v(h(\cdot,s,u)))dv=0.
\eeqnn
For the last two terms on the right hand of (\ref{3.1}), we know from (\ref{1.5}) that
\beqnn
V_t|f|(x)\leq T_t|f|(x)~\mbox{for all}~t\geq 0~\mbox{and}~x\in E.
\eeqnn
Thus by (\ref{2.15}), we get
\beqlb\label{3.3}\nonumber
\ar \ar\left|\exp\left\{\int_0^s\eta(V_u(i\theta_2e^{\lambda_1(t+s)/2} f))du+\int_{M(E)^\circ}\int_0^s
\mathbb{P}_\mu\left(\exp\{i\theta_2e^{\lambda_1(t+s)/2}\<f,X_{u}\>\}-1 \right)du H(d\mu)
\right\}\right|\\\nonumber
\ar \ar\lesssim \exp\left\{\theta_2e^{\lambda_1(t+s)/2}\int_0^s\eta(T_u|f|))du+\theta_2e^{\lambda_1(t+s)/2}\int_{M(E)^\circ}\int_0^s
\mathbb{P}_\mu[\<|f|,X_{u}\>]du H(d\mu)
\right\}\\\nonumber
\ar \ar\lesssim\exp\left\{\theta_2e^{\lambda_1t/2}\left(e^{\lambda_1s/2}\|f\|_2m(t_0)+[e^{(\lambda_1- 2\lambda_{\gamma(f)})s/2}+e^{\lambda_1s/2}]\Gamma(a_{t_0}^{1/2})\right)\right\}\\
\ar \ar\rightarrow 1,~\mbox{as}~t\rightarrow\infty.
\eeqlb
Hence by the dominated convergence theorem, we get
\beqnn
\lim_{t\rightarrow\infty}\kappa(\theta_1,\theta_2,s,t)=\mathbb{P}_\mu\left(\exp\{i\theta_1\widetilde{W}_\infty\}\exp\left\{-\frac{1}{2}\theta_2^2
\<V_s,\phi_1\>_m \widetilde{W}_\infty\right\}\right).
\eeqnn
Since $e^{\lambda_1(t+s)}\<\phi_1,Y_{t+s}\>-e^{\lambda_1t}\<\phi_1,Y_t\>\rightarrow 0$ in probability, as $t\rightarrow\infty$,
we easily get that under $\mathbb{P}_\mu$,
\beqnn
U_3(s,t)\ar:=\ar\left(e^{\lambda_1(t+s)}\<\phi_1,Y_{t+s}\>, e^{\lambda_1(t+s)/2}\<f,Y_{t+s}\>-e^{\lambda_1(t+s)/2}\<T_sf,Y_t\>\right)\\
\ar\stackrel{d}{\rightarrow}\ar(\widetilde{W}_\infty,\sqrt{\widetilde{W}_\infty}G_1(s)),
\eeqnn
as $t\rightarrow\infty$. By (2.15) in \cite{RSZ15},
\beqnn
\lim_{s\rightarrow\infty}V_s(x)=\lim_{s\rightarrow\infty}e^{\lambda_1s}\mathbb{V}ar_{\delta_x}\<f,X_s\>=\sigma_f^2\phi_1(x).
\eeqnn
Thus $\lim_{s\rightarrow\infty}\sigma^2_f(s)=\sigma^2_f$. So $\lim_{s\rightarrow\infty}d(G_1(s),G_1(f))=\sigma^2_f$.

Let $\mathcal{D}(s+t)$ and $\widetilde{\mathcal{D}}(s,t)$ be the distribution of $U_1(s+t)$ and $U_3(s,t)$
respectively, and let $\mathcal{D}(s)$ and $\mathcal{D}$ be the distributions of
$(\widetilde{W}_\infty, \sqrt{\widetilde{W}_\infty}G_1(s))$ and $(\widetilde{W}_\infty, \sqrt{\widetilde{W}_\infty}G_1(f))$
respectively. Then, by a similar argument as in \cite{RSZ14} and using the definition of $\limsup_{t\rightarrow\infty}$,
\beqnn
\limsup_{t\rightarrow\infty}d(\mathcal{D}(t),\mathcal{D})\ar\leq\ar\limsup_{t\rightarrow\infty}d(\mathcal{D}(s+t),\mathcal{D})\\
\ar\leq\ar\limsup_{t\rightarrow\infty}d(\mathcal{D}(s+t),\widetilde{\mathcal{D}}(s,t))+d(\widetilde{\mathcal{D}}(s,t),\mathcal{D}(s))
 +d(\mathcal{D}(s),\mathcal {D})\\
\ar\leq\ar \limsup_{t\rightarrow\infty}(\mathbb{P}_\mu(e^{\lambda_1(t+s)/2}\<T_sf,Y_t\>)^2)^{1/2}+0+d(\mathcal{D}(s),\mathcal {D}).
\eeqnn
Letting $s\rightarrow\infty$, we get
\beqnn
\limsup_{t\rightarrow\infty}d(\mathcal{D}(t),\mathcal{D})\leq\limsup_{s\rightarrow\infty}
\limsup_{t\rightarrow\infty}(\mathbb{P}_\mu(e^{\lambda_1(t+s)/2}\<T_sf,Y_t\>)^2)^{1/2}.
\eeqnn
Therefore, we are left to prove that
\beqlb\label{3.4}
\limsup_{s\rightarrow\infty}\limsup_{t\rightarrow\infty}e^{\lambda_1(t+s)}\mathbb{P}_\mu(\<T_sf,Y_t\>)^2=0.
\eeqlb
By (\ref{2.8}), we have that
\beqnn
\mathbb{P}_\mu(\<T_sf,Y_t\>)^2\ar=\ar\left(\mu(T_{s+t}f)+\int_0^t\Gamma(T_{u+s}f)du\right)^2
+\int_E\int_0^t T_{t-u}[A(T_{u+s}f)^2](x)du\mu(dx)\\ \nonumber
\ar \ar+\int_0^t\int_E\int_0^v T_{v-u}[A(T_{u+s}f)^2](x)du\Gamma(dx)dv\\
\ar \ar+\int_0^t\int_{M(E)^\circ}\nu(T_{u+s}f)^2H(d\nu)du\\
\ar:=\ar B_1+B_2+B_3+B_4.
\eeqnn
By (\ref{2.10}), we get
\beqnn
\lim_{t\rightarrow\infty}e^{\lambda_1(t+s)}B_1=\lim_{t\rightarrow\infty}\left(e^{\lambda_1(t+s)/2}
\left(\mu(T_{s+t}f)+\int_0^t\Gamma(T_{u+s}f)du\right)\right)^2=0
\eeqnn
and we know from the proof of Lemma 3.2 in \cite{RSZ15} that
\beqnn
\lim_{t\rightarrow\infty}e^{\lambda_1(t+s)}B_2=0.
\eeqnn
For $B_3$, by a similar proof as of Lemma 3.2 in \cite{RSZ15}, we get
\beqnn
\ar \ar e^{\lambda_1(t+s)}\int_0^t\int_E\int_0^v T_{v-u}[A(T_{u+s}f)^2](x)du\Gamma(dx)dv\\
\ar \ar= \int_0^te^{\lambda_1(t-v)}e^{\lambda_1(v+s)}\int_E\int_0^v T_{v-u}[A(T_{u+s}f)^2](x)du\Gamma(dx)dv\\
\ar \ar\lesssim \int_0^te^{\lambda_1(t-v)}e^{(\lambda_1-2\lambda_{\gamma(f)})s}dv \Gamma(a_{t_0}^{1/2}).
\eeqnn
As for the last term $B_4$, we have by (\ref{2.6}) that, for $s>t_0$,
\beqnn
\ar \ar e^{\lambda_1(t+s)}\int_0^t\int_{M(E)^\circ}\nu(T_{u+s}f)^2H(d\nu)du\\
\ar \ar \leq e^{(\lambda_1-2\lambda_{\gamma(f)})s}\int_0^te^{(\lambda_1-2\lambda_{\gamma(f)})u}
e^{\lambda_1(t-u)}\int_{M(E)^\circ}\nu(e^{\lambda_{\gamma(f)}(s+u)}T_{s+u}f)^2duH(d\nu)\\
\ar \ar\lesssim e^{(\lambda_1-2\lambda_{\gamma(f)})s}\int_{M(E)^\circ}\nu(a_{t_0}^{1/2})^2H(d\nu).
\eeqnn
Combining the above estimates, we will get the desired result by first letting $t\rightarrow\infty$ and then
letting $s\rightarrow\infty$. \qed

\begin{lemma}\label{l3.4}
If $f\in \mathcal{C}_s$ and $h\in \mathcal{C}_c$. Then
\beqnn
\left(e^{\lambda_1t}\<\phi_1, Y_t\>, t^{-1/2}e^{\lambda_1t/2}\<h,Y_t\>, e^{\lambda_1t/2}\<f,Y_t\> \right)\stackrel{d}\rightarrow
\left(\widetilde{W}_\infty,  G_2(h)\sqrt{\widetilde{W}_\infty}, G_1(f)\sqrt{\widetilde{W}_\infty}\right),
\eeqnn
where $G_1(f)\sim \mathcal{N}(0,\sigma^2_f)$ and $G_2(h)\sim \mathcal{N}(0,\rho^2_h)$. Moreover,
$\widetilde{W}_\infty$, $G_1(f)$ and $G_2(h)$ are independent.
\end{lemma}

\noindent{\it Proof.} We will use
the idea suggested in \cite{RSZ15} with some modifications to prove the result.
In the proof, we always assume $t>3t_0$. We define an $\mathbb{R}^3$-valued random variable by
\beqnn
U_1(t):=\left(e^{\lambda_1t}\<\phi_1,Y_t\>,t^{-1/2}e^{\lambda_1t/2}\<h,Y_t\>, e^{\lambda_1t/2}\<f,Y_t\>\right).
\eeqnn
Let $n>2$ and write
\beqnn
U_1(nt)=\left(e^{\lambda_1nt}\<\phi_1,Y_{nt}\>,(nt)^{-1/2}e^{\lambda_1nt/2}\<h,Y_{nt}\>, e^{\lambda_1nt/2}\<f,Y_{nt}\>\right).
\eeqnn
To consider the limit of $U_1(t)$ as $t\rightarrow\infty$, it is equivalent to consider the limit
of $U_1(nt)$ for any $n>2$. The main idea is as follows. For $t>t_0$, $n>2$,
\beqlb \nonumber
U_1(nt)\ar=\ar\left(e^{\lambda_1nt}\<\phi_1,Y_{nt}\>, \frac{e^{\lambda_1nt/2}(\<h,Y_{nt}\>-\<T_{(n-1)t}h, Y_t\>)}{(nt)^{1/2}},\right.\\\nonumber
\ar \ar\left.e^{\lambda_1nt/2}(\<f,Y_{nt}\>-\<T_{(n-1)t}f, Y_t\>)\right)\\
\ar \ar+\left(0, (nt)^{-1/2}e^{\lambda_1nt/2}\<T_{(n-1)t}h, Y_t\>, e^{\lambda_1nt/2}\<T_{(n-1)t}f, Y_t\>\right).
\eeqlb
We will prove that the second term on the right hand has no contribution to the double limit, first as $t\rightarrow\infty$
and then $n\rightarrow\infty$. The double limit of the first term is equal to another $\mathbb{R}^2$-valued
random variable $U_2(n,t)$ where
\beqnn
U_2(n,t)\ar:=\ar\left(e^{\lambda_1t}\<\phi_1,Y_t\>, \frac{e^{\lambda_1nt/2}(\<h,Y_{nt}\>-\<T_{(n-1)t}h, Y_t\>)}{((n-1)t)^{1/2}},\right.\\
\ar \ar~~~\left.e^{\lambda_1nt/2}(\<f,Y_{nt}\>-\<T_{(n-1)t}f, Y_t\>)\right).
\eeqnn
We claim that
\beqnn
U_2(n,t)\stackrel{d}{\rightarrow}(\widetilde{W}_\infty, \sqrt{\widetilde{W}_\infty}G_2(h),  \sqrt{\widetilde{W}_\infty}G_1(f)),
~~\mbox{as}~t\rightarrow\infty.
\eeqnn
Denote the characteristic function of $U_2(n,t)$ under $\mathbb{P}_\mu$ by $\kappa_1(\theta_1,\theta_2,\theta_3,n,t)$.
Define
\beqnn
Z_1(t,\theta_2):=\theta_2t^{-1/2}e^{\lambda_1t/2}\<h,X_t\>,~~
Z_2(t,\theta_3):=\theta_3e^{\lambda_1t/2}\<f,X_t\>,~t>0,
\eeqnn
and
\beqnn
Z_t(\theta_2,\theta_3):=Z_1(t,\theta_2)+Z_2(t,\theta_3).
\eeqnn
We define the corresponding random variables on $\mathbb{D}$ by $\widetilde{Z}_1(t,\theta_2)$, $\widetilde{Z}_2(t,\theta_3)$
and $\widetilde{Z}_t(\theta_2,\theta_3)$. Using an argument similar to that leading to (\ref{3.1}), we get
\beqlb\label{3.5}\nonumber
\kappa_1(\theta_1,\theta_2,\theta_3,n,t)
\ar=\ar\mathbb{P}_\mu\left(\exp\left\{i\theta_1e^{\lambda_1t}\<\phi_1,Y_t\>+\int_E\int_{\mathbb{D}}
\left(\exp\{ie^{\lambda_1t/2}\widetilde{Z}_{(n-1)t}(\theta_2,\theta_3)(w)\}\right.\right.\right.\\\nonumber
\ar \ar\left.\left.\left.
-1-ie^{\lambda_1t/2}\widetilde{Z}_{(n-1)t}(\theta_2,\theta_3)(w)\right)\mathbb{N}_x(dw)Y_t(dx)\right.\right.\\\nonumber
\ar \ar\left.\left.+\int_t^{nt}\int_{M(E)^\circ}
\mathbb{P}_\mu\left(\exp\{ie^{\lambda_1u/2}Z'_{nt-u}(\theta_2,\theta_3)\}-1\right)H(d\mu)du\right.\right.\\
\ar \ar\left.\left.+\int_0^{(n-1)t}\eta(V_u(i\theta_2[(n-1)t]^{-1/2}e^{\lambda_1nt/2}h+i\theta_3e^{\lambda_1nt/2}f))du\right\}\right),
\eeqlb
where $Z'_{nt-u}(\theta_2,\theta_3)(w)=(\frac{nt-u}{(n-1)t})^{1/2}Z_{1}(nt-u,\theta_2)(w)+Z_{2}(nt-u,\theta_3)(w)$.
Define
\beqnn
R'_t(\theta,x)\ar:=\ar\int_\mathbb{D}\left(\exp\{i\theta \widetilde{Z}_t(\theta_2,\theta_3)(w)\}-1- i\theta \widetilde{Z}_t(\theta_2,\theta_3)(w)\right.\\
\ar \ar\left.+\frac{1}{2}\theta^2(\widetilde{Z}_t(\theta_2,\theta_3)(w))^2\right)\mathbb{N}_x(dw)
\eeqnn
and
\beqnn
J(n,t,x)\ar:=\ar\int_E\int_\mathbb{D}
\left(\exp\{ie^{\lambda_1t/2}\widetilde{Z}_{(n-1)t}(\theta_2,\theta_3)(w)\}\right.\\
\ar \ar\left.-1-ie^{\lambda_1t/2}\widetilde{Z}_{(n-1)t}(\theta_2,\theta_3)(w) \right)\mathbb{N}_x(dw)Y_t(dx).
\eeqnn
Then
\beqnn
J(n,t,x)=-\frac{1}{2}e^{\lambda_1 t}\mathbb{N}_x(\widetilde{Z}_{(n-1)t}(\theta_2,\theta_3))^2+R'_{(n-1)t}(e^{\lambda_1 t/2},x)
\eeqnn
and
\beqnn
\kappa_1(\theta_1,\theta_2,\theta_3,n,t)=\mathbb{P}_\mu(\exp\{i\theta_1e^{\lambda_1t}\<\phi_1,Y_t\>+\<J(n,t,\cdot),Y_t\>\}).
\eeqnn
Let $V_t^n(x):=\mathbb{N}_x(\widetilde{Z}_{(n-1)t}(\theta_2,\theta_3))^2$. Then
\beqnn
\<J(n,t,\cdot),Y_t\>\ar=\ar-\frac{1}{2}e^{\lambda_1 t}\<V_t^n,Y_t\>+\<R'_{(n-1)t}(e^{\lambda_1 t},\cdot),Y_t\>\\
\ar:=\ar J_1(n,t)+J_2(n,t).
\eeqnn
We first consider $J_1(n,t)$.
It follows from the estimate (3.26) in \cite{RSZ15} that
\beqnn
\lim_{t\rightarrow\infty}J_1(n,t)=\lim_{t\rightarrow\infty}-\frac{1}{2}e^{\lambda_1 t}
(\theta_2^2\rho_h^2+\theta_3^2\sigma_f^2)\<\phi_1,Y_t\>=-\frac{1}{2}
(\theta_2^2\rho_h^2+\theta_3^2\sigma_f^2)\widetilde{W}_\infty.
\eeqnn
For $J_2(n,t)$, by an similar argument as in (\ref{3.2}), we have
\beqnn
J_2(n,t)\rightarrow 0 ~~\mbox{in probability, as}~~t\rightarrow\infty.
\eeqnn
We can use the same method as in the proof of (\ref{3.3}) that the last two terms on the right hand side of (\ref{3.5})
tend to
$1$ as $t$ goes to infinity. Hence, combining the above calculations, we get that
\beqnn
\lim_{t\rightarrow\infty}\kappa_1(\theta_1,\theta_2,\theta_3,n,t)
=\mathbb{P}_\mu\left[\exp\{i\theta_1\widetilde{W}_\infty\}\exp\left\{-\frac{1}{2}
(\theta_2^2\rho_h^2+\theta_3^2\sigma_f^2)\widetilde{W}_\infty\right\}\right].
\eeqnn
Similarly as in the proof of Lemma \ref{l3.3}, in order to get the desired result, we now need only to show that
\beqnn
\ar \ar\lim_{n\rightarrow\infty}\limsup_{t\rightarrow\infty}(nt)^{-1}e^{\lambda_1nt}\mathbb{P}_\mu(\<T_{(n-1)t}h, Y_t\>)^2=0,\\
\ar \ar\lim_{n\rightarrow\infty}\limsup_{t\rightarrow\infty}e^{\lambda_1nt}\mathbb{P}_\mu(\<T_{(n-1)t}f, Y_t\>)^2=0.
\eeqnn
In fact, using (\ref{2.12}) and the fact that $\<T_th,\mu\>=e^{-\lambda_1t/2}\<h,\mu\>$, we have
\beqnn
(nt)^{-1}e^{\lambda_1nt}\mathbb{P}_\mu(\<T_{(n-1)t}h, Y_t\>)^2\ar=\ar(nt)^{-1}e^{\lambda_1t}\mathbb{V}\mbox{ar}_\mu\<h, Y_t\>
+(nt)^{-1}e^{\lambda_1t}(\mathbb{P}_\mu\<h, Y_t\>)^2\\
\ar\lesssim\ar n^{-1}(1+t^{-1}).
\eeqnn
Using the same method as the proof of (\ref{3.4}) with $s=(n-1)t$, and then letting $t\rightarrow\infty$, we get
\beqnn
e^{\lambda_1nt}\mathbb{P}_\mu(\<T_{(n-1)t}f, Y_t\>)^2\rightarrow 0.
\eeqnn
The proof is now complete.\qed

Recall that
\beqnn
g(x)=\sum_{k:2\lambda_k<\lambda_1}\sum_{j=1}^{n_k}a_j^k\phi_j^{(k)}(x)~~\mbox{and}
~~I_ug(x)=\sum_{k:2\lambda_k<\lambda_1}\sum_{j=1}^{n_k}e^{\lambda_ku}a_j^k\phi_j^{(k)}(x).
\eeqnn
Note that the sum over $k$ is a sum over a finite number of elements. Define
\beqnn
H_\infty(w):=\sum_{k:2\lambda_k<\lambda_1}\sum_{j=1}^{n_k}a_j^k\widetilde{H}_\infty^{k,j}(w),~w\in\mathbb{D}.
\eeqnn
By Lemma \ref{l3.2}, we have, as $u\rightarrow\infty$
\beqnn
\<I_ug,w_u\>\rightarrow H_\infty,~~\mathbb{N}_x\mbox{-a.e., in}~L^1(\mathbb{N}_x)~\mbox{and in}~L^2(\mathbb{N}_x).
\eeqnn
Since $\mathbb{N}_x\<I_ug,w_u\>=\mathbb{P}_{\delta_x}\<I_ug,X_u\>=g(x)$, we get
\beqnn
\mathbb{N}_x(H_\infty)=g(x)
\eeqnn
and by (3.37) and (3.38) in \cite{RSZ15},
\beqlb\label{3.7}\nonumber
\mathbb{N}_x(H_\infty)^2\ar=\ar\int_0^\infty T_s\left[A\left(\sum_{k:2\lambda_k<\lambda_1}\sum_{j=1}^{n_k}e^{\lambda_ks}a_j^k\phi_j^{(k)}\right)^2 \right](x) ds\\
\ar\lesssim\ar (a_{t_0}(x))^{1/2}\in L^2(E,m)\cap L^4(E,m).
\eeqlb

\medskip

\noindent{\bf Proof of Theorem 3.2.}
Consider an $\mathbb{R}^4$-valued random variable $U_4(t)$ defined by:
\beqnn
U_4(t)\ar:=\ar\left(e^{\lambda_1t}\<\phi_1,Y_t\>, e^{\lambda_1t/2}\left(\<g,Y_t\>-
\sum_{k:2\lambda_k<\lambda_1}\sum_{j=1}^{n_k}e^{-\lambda_kt}a_j^kW_\infty^{k,j}\right),\right.\\
\ar \ar\left.t^{-1/2}e^{\lambda_1t/2}\<h,Y_t\>, e^{\lambda_1t/2}\<f,Y_t\>\right).
\eeqnn
To get the conclusion of Theorem 3.2, it suffice to show that, under $\mathbb{P}_\mu$,
\beqnn
U_4(t)\stackrel{d}{\rightarrow}\left(\widetilde{W}_\infty, \sqrt{\widetilde{W}_\infty}G_3(g),
\sqrt{\widetilde{W}_\infty}G_2(h), \sqrt{\widetilde{W}_\infty}G_1(f)\right),
\eeqnn
where $\widetilde{W}_\infty$, $G_3(g)$, $G_2(h)$ and $G_1(f)$ are independent. Note that,
by Lemma \ref{l3.1},
\beqnn
\lim_{u\rightarrow\infty}\<I_ug,Y_{t+u}\>=\sum_{k:2\lambda_k<\lambda_1}\sum_{j=1}^{n_k}a_j^kW_\infty^{k,j}(w),
~\mathbb{P}_\mu\mbox{-a.s.}
\eeqnn
Denote the
characteristic function of $U_4(t)$ under $\mathbb{P}_\mu$ by $\kappa_2(\theta_1,\theta_2,\theta_3,\theta_4,t)$.
Then we have
\beqlb\label{3.10} \nonumber
\kappa_2(\theta_1,\theta_2,\theta_3,\theta_4,t)\ar=\ar\lim_{u\rightarrow\infty}\mathbb{P}_\mu\left(\exp\left\{i\theta_1e^{\lambda_1t}\<\phi_1,Y_t\>+
i\theta_2e^{\lambda_1t/2}(\<g,Y_t\>-\<I_ug,Y_{t+u}\>)\right.\right.\\ \nonumber
\ar \ar +\left.\left.
i\theta_3t^{-1/2}e^{\lambda_1t/2}\<h,Y_t\>+i\theta_4e^{\lambda_1t/2}\<f,Y_{t}\>\right\}\right)\\ \nonumber
\ar=\ar\lim_{u\rightarrow\infty}\mathbb{P}_\mu\left(\exp\left\{i\theta_1e^{\lambda_1t}\<\phi_1,Y_t\>+
i\theta_3t^{-1/2}e^{\lambda_1t/2}\<h,Y_t\>\right.\right.\\\nonumber
\ar \ar\left.\left.+i\theta_4e^{\lambda_1t/2}\<f,Y_{t}\>+\<J_u(t,\cdot),Y_t\>\right.\right.\\\nonumber
\ar \ar\left.\left.+\int_0^u\int_{M(E)^\circ}
Q_\mu\left(\exp\{-i\theta_2e^{\lambda_1t/2}\<I_ug,X_{u-s}\>\}-1\right)dsH(d\mu)\right.\right.\\
\ar \ar\left.\left.+\int_0^u\eta(V_s(-i\theta_2e^{\lambda_1t/2}I_ug))ds\right\}\right),
\eeqlb
where
\beqnn
J_u(t,x)=\int_\mathbb{D}\left(\exp\{-i\theta_2e^{\lambda_1t/2}\<I_ug,w_u\>\}-1+i\theta_2e^{\lambda_1t/2}\<I_ug,w_u\> \right)\mathbb{N}_x(dw).
\eeqnn
By (3.42) in \cite{RSZ15}, we have
\beqnn
\lim_{u\rightarrow\infty}J_u(t,x)=\mathbb{N}_x
\left(\exp\{-i\theta_2e^{\lambda_1t/2}H_\infty\}-1+i\theta_2e^{\lambda_1t/2}H_\infty\right):=J(t,x)
\eeqnn
and
\beqnn
\lim_{u\rightarrow\infty}\<J_u(t),Y_t\>=\<J(t,\cdot),Y_t\>,~\mathbb{P}_\mu\mbox{-a.s.}
\eeqnn
For the last two terms on the right hand side of (\ref{3.10}), we have
\beqnn
\ar \ar\left|\exp\left\{\int_0^u\eta(V_s(-i\theta_2e^{\lambda_1t/2} I_u g))ds+\int_{M(E)^\circ}\int_0^u
\mathbb{P}_\mu\left(\exp\{-i\theta_2e^{\lambda_1t/2}\<I_ug,X_{s}\>\}-1 \right)ds H(d\mu)
\right\}\right|\\
\ar \ar\lesssim \exp\left\{\theta_2e^{\lambda_1t/2}\int_0^u\eta(T_s|I_ug|))ds+\theta_2e^{\lambda_1t/2}\int_{M(E)^\circ}\int_0^u
\mathbb{P}_\mu[\<|I_ug|,X_{s}\>]ds H(d\mu)
\right\}\\
\ar \ar=\exp\left\{\theta_2e^{\lambda_1t/2}\int_0^u\Gamma(T_s|I_ug|)ds\right\}\\
\ar \ar\leq\exp\left\{\theta_2e^{\lambda_1t/2}\sum_{k:2\lambda_k<\lambda_1}\sum_{j=1}^{n_k}|a_j^k|e^{\lambda_ku}\int_0^u\Gamma(T_s|\phi_j^{(k)}|)ds\right\}\\
\ar \ar\rightarrow 1
\eeqnn
by first letting $u\rightarrow\infty$, then letting $t\rightarrow\infty$ and using (\ref{2.15}). Consequently, by the dominated convergence theorem, we obtain
\beqnn
\kappa_2(\theta_1,\theta_2,\theta_3,\theta_4,t)\ar=\ar\mathbb{P}_\mu\left(\exp\left\{i\theta_1e^{\lambda_1t}\<\phi_1,Y_t\>+
i\theta_3t^{-1/2}e^{\lambda_1t/2}\<h,Y_t\>\right.\right.\\
\ar \ar \left.\left.
+i\theta_4e^{\lambda_1t/2}\<f,Y_{t}\>+\<J(t,\cdot),Y_t\>\right\}\right).
\eeqnn
Let
\beqnn
R_s(\theta,x)=\mathbb{N}_x\left(\exp\{i\theta H_\infty\}-1- i\theta H_\infty+\frac{1}{2}\theta^2H^2_\infty\right).
\eeqnn
Then,
\beqnn
\<J(t),Y_t\>=-\frac{1}{2}\theta_2^2e^{\lambda_1t}\<V,Y_{t}\>+\<R(-e^{\lambda_1t/2}\theta_2,\cdot),Y_t\>,
\eeqnn
where $V(x):=\mathbb{N}_x(H_\infty)^2$. By the results in \cite{RSZ15} and similar arguments as (\ref{3.2}), we have
\beqlb\label{3.39}
\lim_{t\rightarrow\infty}\<R(-e^{\lambda_1t/2}\theta_2,\cdot),Y_t\>=0~~\mbox{in probability.}
\eeqlb
Since $V\in L^2(E,m)\cap L^4(E,m)$, we have by Remark 3.2 that
\beqlb\label{3.40}
\lim_{t\rightarrow\infty}e^{\lambda_1t}\<V,Y_{t}\>=\<V,\phi_1\>_m\widetilde{W}_\infty~~\mbox{in probability.}
\eeqlb
Therefore, combing (\ref{3.39}) and (\ref{3.40}), we get
\beqlb\label{3.41}
\lim_{t\rightarrow\infty}\exp\{\<J(t,\cdot),Y_{t}\>\}=\exp\left\{-\frac{1}{2}\theta_2^2\<V,\phi_1\>_m\widetilde{W}_\infty\right\}~~\mbox{in probability.}
\eeqlb
Recall that $\lim_{t\rightarrow\infty}e^{\lambda_1t}\<\phi_1,Y_t\>=\widetilde{W}_\infty$, $\mathbb{P}_\mu$-a.s.
Thus by (\ref{3.41}) and the fact that $|\exp\{\<J(t,\cdot),Y_{t}\>\}|\leq 1$, we get that as $t\rightarrow\infty$,
\beqnn
\ar \ar\left|\mathbb{P}_\mu\left(\exp\left\{\left(i\theta_1-\frac{1}{2}\theta_2^2\<V,\phi_1\>_m\right)e^{\lambda_1t}\<\phi_1,Y_t\>+
i\theta_3t^{-1/2}e^{\lambda_1t/2}\<h,Y_t\>\right.\right.\right.\\
\ar \ar~~\left.\left.\left.
+i\theta_4e^{\lambda_1t/2}\<f,Y_{t}\>\right\}\right)-\kappa_2(\theta_1,\theta_2,\theta_3,\theta_4,t)\right|\\
\ar \ar\leq \left|\mathbb{P}_\mu\exp\{\<J(t,\cdot),Y_{t}\>\}-\exp\left\{-\frac{1}{2}\theta_2^2\<V,\phi_1\>_me^{\lambda_1t}\<\phi_1,Y_t\>\right\}\right|
\rightarrow 0.
\eeqnn
Consequently, by Lemma \ref{l3.4}, we have
\beqnn
\ar \ar\lim_{t\rightarrow\infty}\mathbb{P}_\mu\left(\exp\left\{\left(i\theta_1-\frac{1}{2}\theta_2^2\<V,\phi_1\>_m\right)e^{\lambda_1t}\<\phi_1,Y_t\>+
i\theta_3t^{-1/2}e^{\lambda_1t/2}\<h,Y_t\>\right.\right.\\
\ar \ar~~ \left.\left.
+i\theta_4e^{\lambda_1t/2}\<f,Y_{t}\>\right\}\right)\\
\ar \ar=\mathbb{P}_\mu\left[\exp\{i\theta_1\widetilde{W}_\infty\}\exp\left\{-\frac{1}{2}
(\theta_2^2\<V,\phi_1\>_m+\theta_3^2\rho_h^2+\theta_4^2\sigma_f^2)\widetilde{W}_\infty\right\}\right].
\eeqnn
By (\ref{3.7}), we get
\beqnn
\<V,\phi_1\>_m=\int_0^\infty e^{\lambda_1s}\<A(I_sg)^2,\phi_1\>_mds.
\eeqnn
The proof is now complete.\qed

\begin{acknowledgements}The parts of this paper were written while the author visited Concordia.
The author would like to give sincere thanks to
Professor Xiaowen Zhou for his encouragement and helpful
discussions and hospitality at Concordia. The author
also thank an anonymous referee for useful comments.
\end{acknowledgements}

\end{document}